\documentclass[twocolumn]{autart}                         
\newtheorem{assumption}{Assumption}[section]

\newtheorem{theorem}{Theorem}[section]
\newtheorem{lemma}{Lemma}[section]

\newtheorem{property}{Property}[section]
\usepackage{amssymb}
\usepackage[fleqn]{amsmath}
\usepackage{amstext}
\usepackage{comment}
\usepackage[colorlinks=true,linkcolor=blue,citecolor=blue,urlcolor=blue]{hyperref}     
\usepackage{float}
\usepackage{graphicx}
\usepackage{soul}
\usepackage{color}
\usepackage{epsfig,amsfonts,color}
\usepackage{bm}
\usepackage{mathrsfs}
\usepackage{multirow}
\usepackage{cite}
\usepackage{algorithmic}
\usepackage{lineno}
\usepackage{enumerate}
\usepackage{array}
\usepackage{placeins}
\usepackage{caption}
\usepackage{subcaption}

\newcolumntype{L}[1]{>{\raggedright\let\newline\\\arraybackslash\hspace{0pt}}m{#1}}
\newcolumntype{C}[1]{>{\centering\let\newline\\\arraybackslash\hspace{0pt}}m{#1}}
\newcolumntype{R}[1]{>{\raggedleft\let\newline\\\arraybackslash\hspace{0pt}}m{#1}}

\allowdisplaybreaks

\begin{document}
\begin{frontmatter}
\title{Hierarchical MPC Schemes for Periodic Systems using Stochastic Programming\thanksref{footnoteinfo}} % Title, preferably not more 
                                                % than 10 words.

\thanks[footnoteinfo]{Corresponding author V.M. Zavala Tel. +1-608-890-2111.}

\author[uw]{Ranjeet Kumar}\ead{rkumar32@wisc.edu},    % Add the 
\author[jci]{Michael J. Wenzel}\ead{mike.wenzel@jci.com},    % e-mail address 
\author[jci]{Matthew J. Ellis}\ead{matthew.j.ellis@jci.com},    % (ead) as shown
\author[jci]{Mohammad N. ElBsat}\ead{mohammad.elbsat@jci.com},    
\author[jci]{Kirk H. Drees}\ead{kirk.h.drees@jci.com},   
\author[uw]{Victor M. Zavala}\ead{victor.zavala@wisc.edu},               

\address[uw]{Department of Chemical and Biological Engineering, University of Wisconsin-Madison\\
1415 Engineering Dr, Madison, WI 53706}
\address[jci]{Johnson Controls International, 507 E. Michigan St., Milwaukee, WI}  % Please supply                                              
             % full addresses      % here.


\begin{abstract}                          
We show that stochastic programming (SP) provides a framework to design hierarchical model predictive control (MPC) schemes for periodic systems. This is based on the observation that, if the state policy of an infinite-horizon problem is periodic, the problem can be cast as a stochastic program (SP). This reveals that it is possible to update periodic state targets by solving a retroactive optimization problem that progressively accumulates historical data. Moreover, we show that the retroactive problem is a statistical approximation of the SP and thus delivers optimal targets in the long run. Notably, this optimality property can be achieved without the need for data forecasts and cannot be achieved by any known proactive receding horizon scheme. The SP setting also reveals that the retroactive problem can be seen as a high-level hierarchical layer that provides targets to guide a low-level MPC controller that operates over a short period at high resolution. We derive a retroactive scheme tailored to linear systems by using cutting plane techniques and suggest strategies to handle nonlinear systems.  
\end{abstract}
\end{frontmatter}

\section{Introduction}

A well-known challenge arising in model predictive control (MPC) is the computational tractability associated with the length of the planning horizon and with the time resolution of the state and control policies \cite{mpctd}. These tractability issues are often encountered in energy system applications that exhibit phenomena and disturbances emanating at multiple timescales. For instance, in energy systems, long horizons are often required to respond to low-frequency (e.g., seasonal) supply/demand variations and peak electricity costs (e.g., demand charges) while fine time resolutions are needed to modulate high-frequency variations (e.g., from wind/solar supply) and to participate in real-time markets \cite{Dowling2017,Braun1990}. Tractability issues are currently handled by using ad-hoc receding horizon (RH) approximations, which are practical but do not provide optimality guarantees \cite{jci1}. 

Hierarchical MPC schemes \cite{scatto1,scatto2} have been recently proposed to handle multiple scales and achieve stability. These schemes, however, do not provide optimality guarantees in the sense that the computed policies match those of the long-horizon problem of interest. The hierarchical scheme proposed in \cite{zavala2016new} uses adjoint information obtained from a long-term and coarse controller to guide a short-term controller operating at fine time resolutions. Computational experiments are provided to demonstrate that this approach can achieve optimality but no guarantees are given. Moreover, such an approach requires smoothness and continuity of the adjoint profiles, which is not guaranteed in general applications.

The hierarchical scheme proposed in this work relies on the observation that, if the optimal policy of an infinite horizon problem is periodic (or can be approximated with a periodic policy), the problem can be cast as a stochastic programming (SP) problem. Periodicity is a property that is commonly observed in systems driven by exogenous factors with strong periodic components (e.g., energy demands and prices) \cite{bieglerperiodic,risbeck2015cost}. Under the SP abstraction, the periodic states are interpreted as design variables and operational policies over the periods are interpreted as recourse variables. We have recently observed that the SP representation provides a mechanism to construct hierarchical MPC schemes in which a {\em long-term} (supervisory) MPC controller provides periodic targets to guide a {\emph short-term MPC} controller \cite{zavalaacc}. Under nominal conditions with {\em perfect forecasts}, we have shown that the hierarchical scheme delivers an optimal policy for the infinite horizon problem. For the more practical case of {\em imperfect forecasts}, the hierarchical scheme needs to re-compute periodic targets. While this can certainly be done using an RH scheme (e.g., computes targets by anticipating multiple future periods), such an approach would not provide optimality guarantees. In fact, to the best of our knowledge, no RH scheme currently exists that can provide optimal policies in the presence of imperfect forecast information. 

The contribution of this work is the observation that, under a periodic setting, one can derive {\em retroactive} schemes that accumulate real historical data to asymptotically deliver optimal targets. We argue that this retroactive design principle offers a key advantage over proactive RH schemes (which rely only on approximate forecasts). The targets obtained with a retroactive scheme are used to guide a low-level controller operating at fine time resolutions within the periods. In the case of linear systems, one can derive a specialized retroactive scheme by using \emph{incremental cutting-plane (CP) techniques} \cite{suvrajeet}. The SP setting also reveals strategies to construct retroactive schemes for nonlinear systems. We demonstrate the concepts using a battery application. 

The paper is structured as follows. In Section \ref{sec:model}, we provide basic definitions. In Section \ref{sec:cuttingplane}, we introduce the concept of retroactive optimization, derive an incremental CP scheme for linear systems, discuss implementation details, and discuss extensions for nonlinear systems. Computational experiments are presented in Section \ref{sec:experiments}. 

%%%%%%%%%%%%%%%%%%%%%%%%%%%%%%%%%%%%%%%%%%
%%%%%%%%%%%%%%%%%%%%%%%%%%%%%%%%%%%%%%%%%%
\section{Basic Definitions and Setting}\label{sec:model}
Consider a long-horizon problem $\mathbf{O}_m$:
\begin{subequations}\label{eq:periodicplanning} \vspace{-0.05in}
\begin{align}
&\nonumber \min_{u_{\xi},x_{\xi},\eta,x_0}  \; \frac{1}{m}\sum_{\xi\in\Xi}\sum_{t\in\mathcal{T}} \varphi_1(x_{\xi,t},u_{\xi,t},d_{\xi,t})+ \eta\\
&\textrm{s.t.} \; \varphi_2(x_{\xi,t},u_{\xi,t},d_{\xi,t})\leq \eta,\; \xi\in\Xi,t\in\mathcal{T}\label{eq:peak2}\\
&\; x_{\xi,t+1}=f(x_{\xi,t},u_{\xi,t},d_{\xi,t}),\; \xi \in \Xi, t\in\mathcal{\bar{T}}\\
&\; x_{\xi+1,0}=x_{\xi,N},\; \xi \in \bar{\Xi}\label{eq:const}\\ 
&\;x_{0,0}= {x}_0\label{eq:per1}\\ 
&\; x_{\xi,t}\in\mathcal{X}, u_{\xi,t}\in\mathcal{U}.  
\end{align}
\end{subequations}
The time horizon is partitioned into {\em periods} $\Xi:=\{1,\dots,m\}$ with intra-period times $\mathcal{T}:=\{0,\dots,n\}$. For convenience, we define the sets $\bar{\mathcal{T}}:=\mathcal{T}\setminus\{n\}$ and $\bar{{\Xi}}:=\Xi\setminus\{m\}$. The controls, states, and data trajectories at period $\xi\in\Xi$ and time $t\in\mathcal{T}$ are denoted as $u_{\xi,t} \in \mathbb{R}^{n_u}, x_{\xi,t}  \in \mathbb{R}^{n_x}$, and $d_{\xi,t} \in \mathbb{R}^{n_d}$, respectively. We define the notation $u_\xi:=\{u_{\xi,t}\}_{t\in\mathcal{T}}$, $x_\xi:=\{x_{\xi,t}\}_{t\in\mathcal{T}}$, and $d_\xi:=\{d_{\xi,t}\}_{t\in\mathcal{T}}$ to denote the inner period trajectories. The problem data can be interpreted as exogenous disturbances or factors affecting the system (e.g., market prices, demands, weather, model errors).

The mapping $\varphi_1: \mathbb{R}^{n_u} \times \mathbb{R}^{n_x} \times \mathbb{R}^{n_d} \rightarrow \mathbb{R}$ is a time-additive cost function and the mapping $\varphi_2: \mathbb{R}^{n_u} \times \mathbb{R}^{n_x} \times \mathbb{R}^{n_d} \rightarrow \mathbb{R}$ is a time-max (peak) cost function. Minimizing the variable $\eta\in \mathbb{R}$ subject to the constraints \eqref{eq:peak2} is equivalent to minimize the peak cost $\max_{\xi\in{\Xi}}\max_{t\in\mathcal{T}}\varphi_2(x_{\xi,t},u_{\xi,t},d_{\xi,t})$. The mapping $f: \mathbb{R}^{n_u} \times \mathbb{R}^{n_x} \times \mathbb{R}^{n_d} \rightarrow \mathbb{R}^{n_x}$ describes the system dynamics and $\mathcal{X}$ and $\mathcal{U}$ are non-empty feasible sets for states and controls, respectively. For reasons that will become apparent, the initial state $x_0\in\mathbb{R}^{n_x}$ in the above formulation is treated as a decision variable. The infinite-horizon problem $\mathbf{O}_\infty$ is obtained by setting $\lim_{m\to \infty}$. If the data $d_\xi,\xi\in \Xi$ is known, the solution of $\mathbf{O}_\infty$ provides a policy with the best possible performance. 

Problem $\mathbf{O}_m$ suffers from tractability issues for large $m$ (long horizons) and/or for large $n$ (fine resolutions). This tractability issue is often handled by using a proactive RH scheme, which approximates the long-horizon policy by moving forward in time and by planning over short time horizons (e.g., that span multiple periods). Our goal is to derive alternative schemes for the case in which an optimal policy of $\mathbf{O}_\infty$ is {\em periodic} (or when this can be approximated reasonably well using a periodic policy).

\begin{assumption}
An optimal policy of $\mathbf{O}_\infty$ satisfies the periodicity constraints $x_{\xi,0}=x_{\xi,n}$ for all $\xi \in {\Xi}$.
\end{assumption}

To obtain an optimal periodic policy, we enforce the periodicity constraints at every period. These constraints, together with the continuity constraints \eqref{eq:const}, can be expressed as $x_{\xi+1,0}=x_{\xi,0},\, \xi\in\bar{\Xi}$. This set of constraints, in turn, can be reformulated as $x_{\xi,0}=x_{0},\, \xi\in{\Xi}$. These modifications give the periodic problem $\mathbf{P}_m$:
\begin{subequations}\label{eq:periodic1} \vspace{-0.05in}
\begin{align}
&\min_{u_{\xi},x_\xi,\eta,x_0}  \; \frac{1}{m}\sum_{\xi\in\Xi}\sum_{t\in\mathcal{T}} \varphi_1(x_{\xi,t},u_{\xi,t},d_{\xi,t})+ \eta\\
&\textrm{s.t.} \; \varphi_2(x_{\xi,t},u_{\xi,t},d_{\xi,t})\leq \eta, \; \xi\in \Xi,t\in\mathcal{T}\\
&\; x_{\xi,t+1}=f(x_{\xi,t},u_{\xi,t},d_{\xi,t}),\; \xi \in \Xi, t\in\mathcal{\bar{T}}\\
&\; x_{\xi,0}=x_{0},\; \xi \in {\Xi}\label{eq:nonant}\\ 
&\; x_{\xi,N}=x_0,\; \xi \in \Xi\\\  
&\; x_{\xi,t}\in\mathcal{X}, u_{\xi,t}\in\mathcal{U}.  
\end{align}
\end{subequations}
The goal of $\mathbf{P}_m$ is to find a periodic state $x_{0}$, peak cost $\eta$, and intra-period policies $u_{\xi},x_{\xi}\; \xi \in \Xi$ that minimize the time-additive and peak costs. We note that, by construction, $x_{\xi+1,0}=x_{\xi,0}=x_0$ holds for all $\xi\in \Xi$ at any solution of $\mathbf{P}_m$. If we have perfect knowledge of the data $d_\xi,\xi\in \Xi$, the infinite-horizon problem $\mathbf{P}_\infty$ identifies a periodic policy with the best possible performance. We denote the solution of $\mathbf{P}_\infty$ as $w^*=(x_0^*,\eta^*)$, and $u^*_\xi,x^*_\xi$. 

Clearly, the feasible region of $\mathbf{P}_m$ is smaller than that of $\mathbf{O}_m$ (the latter does not enforce periodicity). Consequently, the performance of $\mathbf{P}_m$ is expected to be inferior to that of $\mathbf{O}_m$ for all $m$. In some applications, however, the deterioration of performance might be acceptable, and thus a periodic policy will suffice. However, $\mathbf{P}_m$ also becomes intractable for large $m$ and/or $n$. Here, one can address tractability by using a proactive RH scheme, as proposed in \cite{bieglerperiodic}. We will see that this is not necessary, because periodicity results in a structure that enables the derivation of retroactive schemes that are tractable and that offer optimality guarantees. Notably, we will also discuss why such optimality guarantees cannot be provided by proactive schemes.

The structure of $\mathbf{P}_m$ reveals that the only coupling between periods arises from the variables $x_{0}$ and $\eta$. Consequently, by fixing these variables, we can decompose \eqref{eq:periodic1} into individual {\em period subproblems} of the form:
\begin{subequations}\label{eq:period} \vspace{-0.05in}
\begin{align}
\min_{u_{\xi},x_\xi} & \; \sum_{t\in\mathcal{T}} \varphi_1(x_{\xi,t},u_{\xi,t},d_{\xi,t})+ \eta\\
\textrm{s.t.} &\; \varphi_2(x_{\xi,t},u_{\xi,t},d_{\xi,t})\leq \eta,\; t\in\mathcal{T}\label{eq:peak}\\
&\; x_{\xi,t+1}=f(x_{\xi,t},u_{\xi,t},d_{\xi,t}),\; t\in\mathcal{\bar{T}}\\
&\; x_{\xi,0}=x_{0}\\
&\; x_{\xi,N}=x_{0}\\  
&\; x_{\xi,t}\in\mathcal{X}, u_{\xi,t}\in\mathcal{U}.  
\end{align}
\end{subequations}
We define this problem as $\mathbf{S}_\xi$ and we define its optimal cost as $h(w,d_\xi)$. In the following discussion, we use the short-hand notation $h_\xi:=h(w,d_\xi)$. We now make the following key assumption regarding the nature of the data and of the cost $h(w,d_\xi)$. 
\begin{assumption}\label{ass:stoch} Assume that $\{d_\xi\}_{\xi=1}^\infty$ is a sequence of independent and identically distributed (i.i.d.) realizations of a continuous random variable $D$ with associated probability space ($\Omega,\mathcal{F},\mathbb{P}$). Moreover, assume that $\Omega$ is finite ($D$ has finite support) and that the cost function $h:\mathcal{W}\times \Omega\to \mathbb{R}$ is continuous and bounded in its domain.
\end{assumption}
The above assumptions imply that $\mathbf{S}_\xi$ has a feasible solution for all $w\in\mathcal{W}$ and $d_\xi,\,\xi\in\Xi$, and that $\{h_\xi\}_0^\infty$ is an i.i.d. sequence in which $h_\xi$ has a finite second moment and a well-defined mean $\mathbb{E}[w,D]$ for all $w\in\mathcal{W}$. Under these conditions, we can establish the following important property (see page 156 in \cite{shapiro2009lectures}).
\begin{property} \label{prop:lln}Under Assumption \ref{ass:stoch}, the function $\frac{1}{m}\sum_{\xi=1}^m h(w,d_\xi)$ converges pointwise with probability one to $\mathbb{E}[h(w,D)]$ as $m\to \infty$. 
\end{property}
Here, $\mathbb{E}[\cdot]$ is the expectation operator. This property is the so-called weak law of large numbers (LLN) and is key because it reveals that the infinite horizon problem $\mathbf{P}_\infty$ can be interpreted as a {\em stochastic programming (SP) problem} of the form: 
\vspace{-0.1in}
\begin{align}\label{eq:twostage}
\min_{w\in \mathcal{W}} \; & \phi(w):=g(w) + \mathbb{E}\left[h(w,D)\right].
\end{align}
\vspace{-0.0in}
Under this representation, periods are realizations $d_\xi$ of $ D$, $w:=(x_{0},\eta)$ are {\em design} (target) variables, and $x_{\xi},u_{\xi}$ are (recourse) policies associated with realization $d_\xi$. Here, $g:\mathcal{W}\to \mathbb{R}$ is a bounded cost function for the design variables. In the context of $\mathbf{P}_\infty$, we have $g(w)=\eta$ but this function can be generalized to enforce a cost also on $x_0$. Since \eqref{eq:twostage} is equivalent to $\mathbf{P}_\infty$, it delivers optimal targets $w^*$. We define the solution set of \eqref{eq:twostage} as $S\subseteq\mathcal{W}$.

\vspace{-0.1in}
{\bf Remark:} The i.i.d. requirement on $\{d_\xi\}$ ensures that the LLN holds. This requirement can be enforced by defining a sufficiently long period duration (to eliminate autocorrelation between periods). This approach is commonly used in long-term statistical extrapolation and time series analysis \cite{box2015time,ragan2008statistical}. The i.i.d. requirement can also be relaxed by allowing for sequences that have bounded correlation \cite{hu2008convergence}. To give an idea of why this is the case, we provide the following basic result. 
\begin{property} Assume that $\{h_\xi\}_{-\infty}^{\infty}$ is a sequence of identically distributed random variables with finite second moment and mean $\mathbb{E}[h(x,D)]$. Assume also that there exists $0\leq c<\infty$ such that $\sum_{\xi=-\infty}^{\infty}\mathbb{E}[(h_{\xi'}-\mu)(h_{\xi}-\mu)]\leq c$ holds for all $\xi'=1,\dots,\infty$. We have that $\lim_{m\to \infty}\frac{1}{m}\sum_{\xi=1}^mh_\xi=\mathbb{E}[h(x,D)]$. 
\end{property}
{\bf Proof:} Define $S_m:=\sum_{\xi=1}^mh_\xi$, $\bar{h}_m:=S_m/m$, and $\mu:=\mathbb{E}[h(x,D)]$. The variance of $S_m$ (denoted as $\mathbb{V}[S_m]$) is:
\begin{align*}
\mathbb{E}[(S_m-m\,\mu)^2]&=\sum_{\xi'=1}^m\sum_{\xi=1}^m\mathbb{E}[(h_{\xi'}-\mu)(h_{\xi}-\mu)]\\
&\leq\sum_{\xi'=1}^{m}\sum_{\xi=-\infty}^{\infty}\mathbb{E}[(h_{\xi'}-\mu)(h_{\xi}-\mu)]\\
&\leq m\,c.
\end{align*}
We thus have $\mathbb{V}[\bar{h}_m]\leq c/m$ and, from Chebyshev's inequality,
\begin{align*}
\mathbb{P}(|\bar{h}_m-\mu|>\kappa)&\leq \frac{1}{\kappa^2}\mathbb{E}[(\bar{h}_m-\mu)^2]\\
&= \frac{1}{\kappa^2m^2}\mathbb{E}[(S_m-m\mu)^2]\\
&\leq \frac{m\,c}{m^2\kappa^2}=\frac{c}{m\,\kappa^2}.
\end{align*}
We thus have $\lim_{m\to\infty}\bar{h}_m= \mu$ with probability one. $\Box$

\section{Hierarchical Schemes based on Statistical Approximations} \label{sec:cuttingplane}

The SP representation opens a number of interesting directions. In particular, it provides a mechanism to derive {\em hierarchical MPC schemes}. For instance, one can use a statistical approximation of \eqref{eq:twostage} (equivalently of $\mathbf{P}_\infty$) to provide targets $x_{0},\eta$ that guide a short-term MPC controller of the form \eqref{eq:period}. As we discuss next, approximations of $\mathbf{P}_\infty$ can be constructed and solved in a tractable manner by using well-established SP tools \cite{zavalalaird,caroe1999dual,benders}.  

\subsection{Retroactive Optimization}

A key observation that arises from the SP representation is that one can derive {\em retroactive optimization} schemes that accumulate data over time to refine targets. To explain how such an scheme would work, assume that the system is currently at the beginning of period $m+1$ and the data history $\{d_\xi\}_{\xi=1}^{m}$ is known. We use this information to solve the approximation of \eqref{eq:twostage}:
\begin{align}\label{eq:twostageapprox}
\min_{w\in \mathcal{W}} \;  \phantom{ } & \phi_m(w):=g(w) + \frac{1}{m}\sum_{\xi=1}^m h(w,d_\xi).
\end{align}
This problem is equivalent to $\mathbf{P}_m$ and, because $\{d_\xi\}_{\xi=1}^m$ is i.i.d., $\mathbf{P}_m$ is an statistical approximation of $\mathbf{P}_\infty$. In the SP literature problem, \eqref{eq:twostageapprox} is known as a sample average approximation (SAA). We define the solution set of $\mathbf{P}_m$ as $S_m\subseteq\mathcal{W}$. A solution of $\mathbf{P}_m$ is used to update the targets for the next period $w_{m+1}=(x_{0,m+1},\eta_{m+1})$. The solution of $\mathbf{P}_m$ also implicitly computes optimal (retroactive) policies $u_{\xi},x_{\xi},\,\xi=1,\dots,m$ associated with the sequence $\{d_\xi\}_{\xi=1}^m$. These retroactive policies are interpreted as policies that the system would have taken given knowledge of the data.  

Given that the system is at the current state target $x_{0,m}$ (obtained from the solution of $\mathbf{P}_{m-1}$), we use the targets $w_{m+1}$ to guide a short-term MPC controller over period $m+1$. At this point, however, the data $d_{m+1}$ is not known, so we use a forecast $\hat{d}_{m+1}$ to find policies that optimize the system over period $m+1$ while satisfying the targets $w_{m+1}$. This can be interpreted as solving $\mathbf{S}_\xi$ for $w_{m+1}$ and $\hat{d}_{m+1}$. At the beginning of the next period $m+2$, the data $d_{m+1}$ reveals itself and we use this to solve the approximation $\mathbf{P}_{m+1}$ to obtain targets $w_{m+2}$. 

The retroactive scheme is consistent because, from Property \ref{prop:lln}, we know that accumulating data over time will yield an asymptotically exact statistical approximation $\lim_{m\to \infty} \; \phi_m(w)=\phi(w)$ and thus the targets obtained with $\mathbf{P}_m$ will provide a solution to $\mathbf{P}_\infty$ as $m\to\infty$. This asymptotic convergence result is formally stated in the following theorem (see Theorem 5.3 in \cite{shapiro2009lectures}).

\begin{theorem} Suppose that there exists a compact set $C\subset\mathbb{R}^{n_w}$ such that: i) the solution set $S$ of $\mathbf{P}_\infty$ is nonempty and contained in $C$, ii) the function $\phi(w)$ is finite valued and continuous on $C$, iii) $\phi_m(w)$ converges to $\phi(w)$ with probability one as $m\to \infty$ uniformly in $w\in C$, and with probability one for large enough $m$ the solution set $S_m$ of $\mathbf{P}_m$ is nonempty and contained in $C$. Then the distance between the solution sets $\mathbb{D}(S_m,S)$ converges to zero with probability one as $m\to \infty$. 
\end{theorem}
The above result implies that the distance of any solution of $\mathbf{P}_m$ to the solution set of $\mathbf{P}_\infty$ is zero as $m\to \infty$. Statistical approximation results for SPs also indicate that the probability that a solution of $\mathcal{P}_m$ is in the solution set of $\mathcal{P}_\infty$ increases exponentially with $m$ (Theorem 5.16 in \cite{shapiro2009lectures}). In other words, the probability of finding better targets than those obtained with $\mathbf{P}_m$ decays exponentially as information is accumulated over time. 

These asymptotic optimality results provide a key advantage of the retroactive scheme over traditional RH schemes. This is based on a fundamental design difference: the retroactive scheme uses past (but real) data while RH schemes use future (but approximate) data forecasts. Moreover, proactive schemes discard historical data when computing new targets. The fact that historical data is discarded prevents RH schemes from offering asymptotic optimality guarantees. The structure of $\mathbf{P}_m$ can be exploited using decomposition strategies and this enables scalability to large values of $m$. In the following section, we provide a specific scheme for linear systems, and we then discuss extensions to nonlinear systems.  

\subsection{Incremental Cutting Plane Scheme}

We derive a retroactive scheme tailored to linear systems that delivers optimal targets in the long run. This is done by using an {\emph {incremental cutting plane}} (CP) scheme that accumulates data over time. Our approach is an adaptation of the stochastic decomposition scheme proposed in \cite{suvrajeet} to tackle linear SPs. To derive this linear setting, we assume that $g(w)=c_w^Tw$ where $c_w\in\mathbb{R}^{n_w}$ is a cost vector, we assume that the set $\mathcal{W}\subseteq\mathbb{R}^{n_w}$ is polyhedral, and we assume that $\mathbf{S}_\xi$ has the form:
\begin{align} \label{eq:subobj}
h(w,d_\xi) :=  & \min_{y\in \mathbb{R}_+^{n_y}} \; c^T_\xi y  \quad  \textrm{ s.t.} \quad W y = r_\xi - Tw. 
\end{align}
Here, the data realization is given by $d_\xi=(c_\xi,r_\xi)$. We use $y(w,d_\xi)$ to denote the primal solution vector containing the intra-period trajectories $(x_{\xi}, u_{\xi})$ and some additional dummy variables. The intra-period dynamics are captured using the $W,T$ are coefficient matrices. The structure of the recourse problem is used to simplify algebraic manipulations and is done without loss of generality (the results that we present hold provided that the recourse problem is a linear program). The representation of $\mathbf{S_\xi}$ allows us to express its dual form in the following compact form:
\vspace{-0.1in}
\begin{align}\label{eq:subprobkdual}
& \max_{\pi} \; \pi^T(r_\xi - Tw) \quad \textrm{s.t.} \quad W^T \pi \leq c_\xi. 
\end{align}
\vspace{-0.3in}

Here, $\pi(w,d_\xi)$ is the dual solution vector (dual variables of \eqref{eq:subprobkdual}) and we recall that $\pi(w,d_\xi)^T(r_\xi - Tw)\leq h(w,d_\xi)$ holds for $w\in \mathcal{W}$. We assume that the feasible set of the dual subproblem is a non-empty, compact, and convex polyhedral set. As a result, the subproblem has a finite number of dual vertices. Moreover, because the support $\Omega$ is finite, the set of vertices is finite. We use $\mathbb{V}$ to denote the set of all dual vertices. We also note that, by definition, $\pi(w,d_\xi)\in \mathbb{V}$ for all $(w,d_\xi)\in \mathcal{W}\times \Omega$. 

In the linear setting, the cost function of $\mathbf{P}_\infty$ (given by $\phi(w)=c_w^Tw+\mathbb{E}[h(w,D)]$) can be outer-approximated using CPs accumulated over $\xi=1,\dots,m$ as:
{\setlength{\mathindent}{0pt}
\begin{align}
\underline{\phi}_m(w) :=   \max \{\alpha_\xi^m +(c_w+\beta_\xi^m)^Tw \vert \xi=1,\dots,m\}, \label{fkx}
\end{align}}
where the coefficients $\alpha_m^m$, $\beta_m^m$ are selected to match:
\begin{align}\label{eq:coeff}
&\alpha_m^m + (c_w+\beta_m^m)^Tw\nonumber \\
&\qquad= c^T_ww + \frac{1}{m} \sum_{\xi=1}^m{({\pi}_\xi^m)}^T (r_\xi- Tw).
\end{align}
Here, ${\pi}_\xi^m \in \textrm{argmax}\{\pi^T(r_\xi - Tw_m) \vert \pi \in \mathbb{V}_m\}$ for $\xi=1,\dots,m$, where $\mathbb{V}_m\subseteq\mathbb{V}$ is the collection of vertices accumulated up to period $m$. For convenience, we define the function $h_m(w,d_\xi):= \max\{\pi^T(r_\xi - Tw) \vert \pi \in \mathbb{V}_m\}$ and note that $h_m(w,d_\xi)=({\pi_\xi^{m}})^T(r_\xi - Tw)$ holds. Moreover, we note that $h(w,d_\xi)= \max\{\pi^T(r_\xi - Tw) \vert \pi \in \mathbb{V}\}$ holds. Since $\mathbb{V}_m\subseteq\mathbb{V}$, we have that $({\pi_\xi^{m}})^T(r_\xi - Tw)\leq h(w,d_\xi)$. 

The {\em running cost} is given by $\phi_m(w)=c_w^Tw+\frac{1}{m} \sum_{\xi=1}^m h(w,d_\xi)$. We will prove that $\underline{\phi}_m(w)$ underestimates the running cost $\phi_m(w)$ for all $m$ and converges to the infinite-horizon cost $\phi(w)$ as $m\to \infty$. Consequently, at time $m$, we update the targets by solving the {\em master} problem $\mathbf{M}_m$:
\vspace{-0.1in}
\begin{align}\label{eq:masterk}
 \min_{w\in \mathcal{W}} \; \underline{\phi}_m(w).
\end{align}
This problem is a tractable surrogate of $ \min_{w\in \mathcal{W}} \; {\phi}_m(w)$, because it captures the recourse subproblems by using hyperplanes. This becomes particularly important as information is accumulated over time. The solution of $\mathbf{M}_m$ is used to update the targets $w_{m+1}$, which in turn are used to solve the recourse subproblem $S_{\xi+1}$ and with this, obtain a new dual vertex to be stored in $\mathbb{V}_{m+1}$. The CP scheme is summarized as:
\vspace{-0.1in}
\begin{enumerate}
\item Initialize $m \leftarrow 1$, $\mathbb{V}_{m} \leftarrow \emptyset$, and $w_{m}$.
\item At period time $m+1$: 
\item Observe $d_{m}$ and solve $\mathbf{S_m}$ to obtain $\pi(w_{m},d_{m})$. 
\item Update $\mathbb{V}_m\leftarrow \mathbb{V}_{m-1} \cup \{\pi(w_{m},d_{m})\}$. 
\item \begin{enumerate}
			\item{Obtain ${\pi}_\xi^m \in \textrm{argmax}\{\pi^T(r_\xi - Tw_m) \vert \pi \in \mathbb{V}_m\}$} for all $\xi=1,\dots,m$.		
			\item Get $\alpha_m^m$ from \eqref{eq:coeff}.
			\item Update $\alpha_\xi^m\leftarrow \frac{m-1}{m}\alpha_\xi, \beta_\xi^m\leftarrow \frac{m-1}{m}\beta_\xi$ for $\xi =1,\dots,m-1$.			
\end{enumerate}
\item Solve $\mathbf{M}_m$ and obtain updated targets $w_{m+1}$.
\item Shift period time $m\leftarrow m+1$ and return to Step 2.
\end{enumerate}

We now prove that the CP scheme delivers a sequence of targets $\{w_m\}_{m=1}^\infty$ that converges to optimal targets $w^*$ of $\mathbf{P}_\infty$. Our analysis follows along the lines of that presented in \cite{suvrajeet}. 
\begin{theorem}
The CP $\alpha_m^m + (c_w+\beta_m^m )^Tw_m$ generated in period $m$ provides a statistically valid lower bound for $\phi(w)$ for all $w\in \mathcal{W}$. 
\end{theorem}
\emph{Proof.} Because $\mathbb{V}_m \subseteq \mathbb{V}$ we have that,
\begin{align*}
& \max\{\pi^T(r_\xi - Tw_m) \vert \pi \in \mathbb{V}_m\}\\
&\qquad\qquad \leq \max\{\pi^T(r_\xi - Tw_m) \vert \pi \in \mathbb{V}\},
\end{align*}
and $({{\pi}_\xi^m})^{T}(r_\xi - Tw_m) \leq h(w_m,d_\xi),\; \xi=1,\dots,m$.
We thus have
\begin{align*}
\frac{1}{m} \sum_{\xi=1}^m ({{\pi}_\xi^m})^{T}(r_\xi - Tw_m) \leq \frac{1}{m} \sum_{\xi=1}^m h(w_m,d_\xi).
\end{align*}
Furthermore, $\pi^T(r_\xi - Tw) \leq h(w,d_\xi)$ for any $\pi\in \mathbb{V}$ and:
\begin{align*} 
&c^T_ww+\frac{1}{m} \sum_{\xi=1}^m ({{\pi}_\xi^m})^{T}(r_\xi - Tw) \nonumber\\
&\qquad \leq 
c^T_ww+\frac{1}{m} \sum_{\xi=1}^m h(w,d_\xi)=\phi_m(w),\; w \in \mathcal{W}.  
\end{align*}
The result follows from \eqref{eq:coeff} and by noticing that $\phi_m(w)$ is a statistical approximation of $\phi(w)$. $\Box$

As more observations $d_\xi$ are collected, it is important that all the collected CPs provide a statistically valid lower bound for $\phi(w)$. This is the goal of Step 5b in the CP scheme. In particular, at period $m+p$ with $p>0$:
\begin{align*}
\nonumber \frac{1}{m+p} \sum_{\xi=1}^m ({{\pi}_\xi^m})^{T}(r_\xi - Tw) & \leq \frac{1}{m+p} \sum_{\xi=1}^m h(w,d_\xi) \\ 
& \leq \frac{1}{m+p} \sum_{\xi=1}^{m+p} h(w,d_\xi)
\end{align*}
Thus, in the $(m+p)^{th}$ period, the CP
\begin{align*}
&c^T_ww + \frac{1}{m+p} \sum_{\xi=1}^{m} ({{\pi}_\xi^m})^{T}(r_\xi - Tw) \\
&\qquad\qquad =  \alpha_m^{(m+p)} + (\beta_m^{(m+p)} + c_w)^Tw
\end{align*}
still provides a statistically valid lower bound for $\phi(w)$. 

We now explore the limiting behavior of $h_m(\cdot)$, which embeds all CP information accumulated over time. 
\begin{lemma}\label{Lem:hkconv}
The sequence $\{h_m(\cdot)\}^\infty_{m=1}$ converges uniformly on $\mathcal{W}$.
\end{lemma}
\emph{Proof.} $\mathbb{V}_m \subseteq \mathbb{V}_{m+1} \subseteq \mathbb{V}$ implies that $h_m(w,d_\xi) \leq h_{m+1}(w,d_\xi) \leq h(w,d_\xi)$ for all $w \in \mathcal{W}$ and $d_\xi$. Since $\{h_m(\cdot)\}^\infty_{k=1}$ increases monotonically and is bounded by the finite function $h(\cdot)$, it follows that $\{h_m(\cdot)\}^\infty_{k=1}$ converges point-wise to some function $\varphi(\cdot)$ satisfying $\varphi(w,d_\xi) \leq h(w,d_\xi)$ for all $w\in\mathcal{W}$, $d_\xi\in\Omega$. Since $\mathbb{V}_m \subseteq \mathbb{V}_{m+1}\subseteq \mathbb{V}$, we have $\bar{\mathbb{V}} = \lim_{m \rightarrow \infty} \mathbb{V}_m \subseteq \mathbb{V}$. Since $\mathbb{V}$ is a finite set, so $\bar{\mathbb{V}}$ is also a finite set, and we have that:
\begin{align*}
\varphi(w,d_\xi) &= \lim_{m \rightarrow \infty} h_m(w,d_\xi) \\
&= \lim_{m \rightarrow \infty} \max\{\pi^T(r_\xi-Tw) \vert \pi \in \mathbb{V}_m\} \\
&= \max\{\pi^T(r_\xi-Tw) \vert \pi \in \bar{\mathbb{V}}\},
\end{align*}
and it follows that $\varphi(\cdot)$ is a continuous function. As $\mathcal{W} \times \Omega$ is a compact space, and $\{h_m(\cdot)\}^\infty_{m=1}$ is a monotonic sequence of continuous functions, it follows that it converges uniformly to $\varphi(\cdot)$ (see Theorem 7.13 in \cite{rudin1964principles}). $\Box$

We now show that the sequence of CPs generated provides support points for $\phi(\cdot)$ in the limit $m\to \infty$. 
\begin{theorem}\label{Thm:xknconv}
Let $\{w_{m_k}\}^\infty_{k=1}$ be an infinite subsequence of $\{w_m\}^\infty_{m=1}$. If $w_{m_k} \rightarrow \bar{w}$ then, with probability one,
\begin{align*}
\lim_{k\to \infty}\frac{1}{m_k} \sum_{\xi=1}^{m_k} ({{\pi}_\xi^{m_k}})^T(r_\xi - Tw_{m_k})=\mathbb{E}[h(\bar{w},D)].
\end{align*}
In addition, every limit of $\{\alpha_{m_k}^{m_k},\beta_{m_k}^{m_k}+c_w\}^\infty_{k=1}$ defines a support of $\phi(w)$ at $\bar{w}$, with probability one.
\end{theorem}
\emph{Proof.} By definition, we have that
{\setlength{\mathindent}{10pt}
\begin{align*}
h_{m_k}(w_{m_k},d_\xi) &= ({\pi^{m_k}_\xi})^T(r_\xi - Tw_{m_k})\\
\frac{1}{m_k} \sum_{\xi=1}^{m_k} h_{m_k}(w_{m_k},d_\xi) &= \frac{1}{m_k} \sum_{\xi=1}^{m_k} ({\pi_\xi^{m_k}})^T(r_\xi - Tw_{m_k}).
\end{align*}}
By Lemma \ref{Lem:hkconv}, there exists $\varphi(\cdot) \leq h(\cdot)$ such that $\{h_m\}^\infty_{m=0}$ converges uniformly to $\varphi(\cdot)$. We thus have:
\begin{align*}
\lim_{k\to\infty}\frac{1}{m_k} \sum_{\xi=1}^{m_k} \left[h_{m_k}(w_{m_k},d_\xi) - \varphi(\bar{w},d_\xi)\right] = 0,
\end{align*}
and $\lim_{k\to\infty}\frac{1}{m_k} \sum_{\xi=1}^{m_k} h(w,d_\xi) = \mathbb{E}[h(w,D)]$ with probability one. It is now sufficient to show that $\varphi(\bar{w},d_\xi)=h(\bar{w},d_\xi)$ with probability one. Let $d_\xi$ be a given realization and suppose that, for every $\delta > 0$, we have $\mathbb{P}\{\vert D-d_\xi \vert < \delta\}$. Then, for every $\delta > 0, {\vert d_{m_k}-d_\xi \vert} < \delta $ infinitely often, with probability one. Because $h(\cdot)$ is a continuous function and $\{h_m(\cdot)\}^\infty_{m=1}$ is uniformly convergent, for every $\epsilon > 0$ there exists a $\delta>0$ and $N<\infty$ such that ${\vert (\bar{w},d_{\xi})-({w},d) \vert} < \delta$ and with:
\begin{align*}
& {\vert h(\bar{w},d_{\xi})-h(\bar{w},d) \vert} < \frac{\epsilon}{3} \\
& {\vert h_{m_k}(\bar{w},d_{\xi})-h_{m_k}(\bar{w},d) \vert} < \frac{\epsilon}{3}, \quad  k \geq N.
\end{align*}
Consequently, because $\lim_{k\to\infty}w_{m_k} = \bar{w}$, we have $\mathbb{P}\{\vert D-d_\xi \vert < \delta\}$ implies that for every $\epsilon > 0$ there exists a further subsequence $\{(w_{m'_k},d_{m'_k})\}^\infty_{k=1}$ and $K<\infty$ such that
\begin{align*}
& {\vert h(\bar{w},d_\xi)-h_{m_k'}(\bar{w},d_{m'_k}) \vert} < \frac{\epsilon}{3} \\
&{\vert h(\bar{w},d_{m'_k})-h(w_{m'_k},d_{m'_k}) \vert} < \frac{\epsilon}{3} \\
&{\vert h_{m'_k}(w_{m'_k},d_{m'_k})-h_{m'_k}(w_{m'_k},d_\xi) \vert} < \frac{\epsilon}{3}
\end{align*}
for all $m'_k \geq K$. By construction, {$h_{m'_k}(w_{m'_k},d_{m'_k})=h(w_{m'_k},d_{m'_k})$}. Thus, for every $\epsilon > 0$, there exists a subsequence $\{w_{m'_k}\}^\infty_{k=1}$ and $K < \infty$ such that
\begin{align*}
&{\vert h(\bar{w},d_{\xi})-h_{m'_k}(w_{m'_k},d_\xi) \vert}  <  \epsilon
\end{align*}
for all $m'_k \geq K$. Consequently, it follows that $\varphi(\bar{w},d_\xi) = h(\bar{w},d_\xi)$. Finally, since $\Omega$ is compact, we have that $\mathbb{P}\{\vert D-d_\xi \vert < \delta\}$ for some $\delta > 0$ and for only finitely many values of $d_\xi$, with probability one. Thus, with probability one, $\varphi(\bar{w},d_\xi) = h(\bar{w},d_\xi)$ for all but a finite number of realizations. Hence, with probability one, 
\begin{align*}
\lim_{k\to\infty}\frac{1}{m_k} \sum_{\xi=1}^{m_k} (\pi_\xi^{m_k})^T(r_\xi - Tw_{m_k}) = \mathbb{E}[h(w,D)].
\end{align*}
Moreover, since $h(w,d_\xi) = \max\{\pi^T(r_\xi-Tw) \vert \pi \in \mathbb{V}\}$ and $\mathbb{V}_m \subset \mathbb{V}$ for all $m$, it follows that 
{\setlength{\mathindent}{0pt}
\begin{align*}
c^T_ww + \frac{1}{m_k} \sum_{\xi=1}^{m_k} h(w,d_\xi) & \geq c^T_ww + \frac{1}{m_k} \sum_{\xi=1}^{m_k} \pi_\xi^{m_k^T}(r_\xi - Tw)\\
&=\alpha^{m_k}_{m_k} + \left(\beta^{m_k}_{m_k}+c_w\right)^Tw.
\end{align*}}
We conclude that, with probability one, $\phi(w)$ is at least as large as any limiting cut that is associated with the subsequence of cuts defined by $\{(\alpha^{m_k}_{m_k}, \beta^{m_k}_{m_k}+c)\}^\infty_{k=1}$. Thus, any limiting cut defines a support of $\phi(w)$ at $\bar{w}$. $\Box$

\begin{theorem} \label{Thm:fxkn}
There exists a subsequence of $\{w_{m_k}\}^\infty_{k=1}$, satisfying $\lim_{k \rightarrow \infty} (\underline{\phi}_{m_k}(w_{m_k})-\underline{\phi}_{m_k-1}(w_{m_k}))=0$.
\end{theorem}
\emph{Proof.} See proof of Theorem 3 in \cite{suvrajeet}. 

We now establish the main convergence result. 
\begin{theorem} \label{Thm:xkn}
There exists a subsequence $\{w_{m_k}\}^\infty_{k=1}$, such that every accumulation point of $\{w_{m_k}\}^\infty_{k=1}$ is an optimal solution $w^*$ of $\mathbf{P}_\infty$, with probability one.
\end{theorem}
\emph{Proof.} Let $w^*$ be an optimal solution of $\mathbf{P}_\infty$. From Theorem \ref{Thm:fxkn}, there exists a subsequence $\{w_{m_k}\}^\infty_{k=1}$ such that $\lim_{k \rightarrow \infty} (\underline{\phi}_{m_k}(w_{m_k})-\underline{\phi}_{m_k-1}(w_{m_k}))=0$. Let $\{w_{m_k}\}_{k\in\mathcal{K}}$ be a further subsequence such that $\lim_{k\in\mathcal{K}} w_{m_k} = \bar{w}$. Compactness of $\mathcal{W}$ ensures that $\bar{w} \in \mathcal{W}$, and thus $\phi(w^*) \leq \phi(\bar{w})$. We know that:
\begin{align} \label{eq:final1}
\underline{\phi}_m(w) \leq c^T_ww + \frac{1}{m} \sum_{\xi=1}^m h(w,d_\xi),\; w\in \mathcal{W},
\end{align}
and thus, 
{\setlength{\mathindent}{0pt}
\begin{align} \label{eq:final2}
\limsup_{m \in \mathcal{M}} \underline{\phi}_m(w^*) \leq c^T_ww^* + \mathbb{E} [h(w^*,D)] = \phi(w^*)
\end{align}}
with probability one for any index set $\mathcal{M}$. Since $w_m$ minimizes $\underline{\phi}_{m-1}(\cdot)$ on $\mathcal{W}$, we have $\underline{\phi}_{m-1}(w_m) \leq \underline{\phi}_{m-1}(w^*)$
for all $m$. From Theorem \ref{Thm:xknconv}, $\lim_{k\in\mathcal{K}} \phi_{m_k}(w_{m_k}) \leq \phi(\bar{w})$ with probability one and, also by definition, $\lim_{k\in\mathcal{K}} \underline{\phi}_{m_k}(w_{m_k}) -  \underline{\phi}_{m_k-1}(w_{m_k}) = 0$. Thus, we have $\lim_{k\in\mathcal{K}} \underline{\phi}_{m_k-1}(w_{m_k}) = \phi(\bar{w})$, with probability one. Combining these results we obtain:
{\setlength{\mathindent}{0pt}
\begin{align*} \label{eq:final4}
\phi(w^*) &\leq \phi(\bar{w}) \nonumber\\
&= \limsup_{k\in\mathcal{K}} \underline{\phi}_{m_k-1}(w_{m_k}) \nonumber\\
&\leq  \underline{\phi}_{m_k-1}(w^*) \leq \phi(w^*),
\end{align*}}
with probability one. We thus have that $\bar{w}=w^*$. $\Box$

\subsection{Short-Term MPC Controller}

The CP scheme is guaranteed to deliver optimal targets as data is accumulated over time. Notably, because the scheme is inherently retroactive, it achieves optimal targets without using any data forecasts. So the question is: How does the proposed scheme accommodates forecast information? From an implementation stand-point, another important question is: What metrics can one use to monitor optimality of the targets? 

\begin{figure}[!htb]
\begin{center}
\includegraphics[width=3.5in]{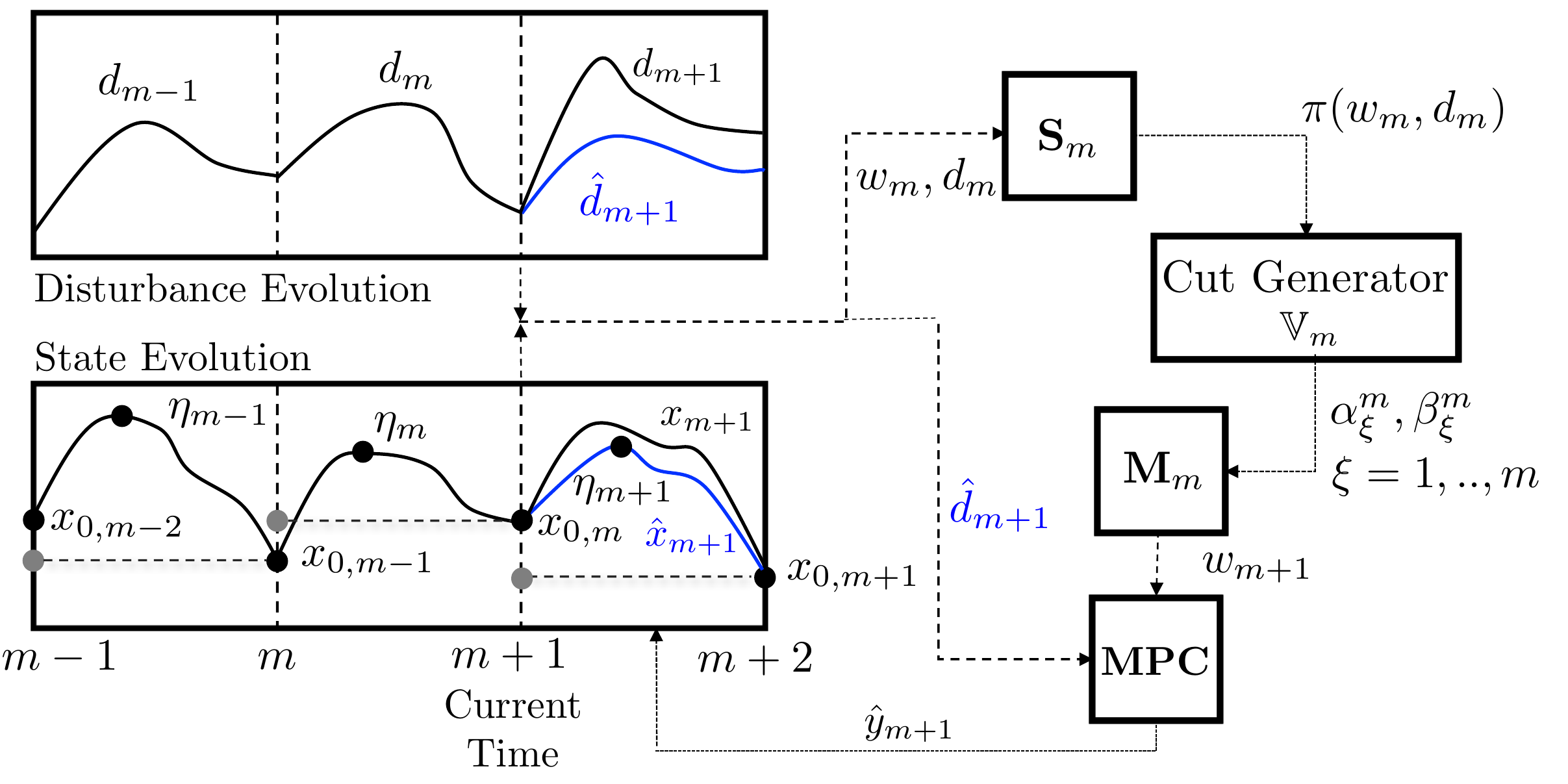}
\caption{Sketch of hierarchical scheme using cutting planes.}\label{fig: summary}
\end{center}
\end{figure}

The proposed scheme offers a couple of mechanisms to embed forecast information. First, initial guesses for the targets $w_m=(x_{0,m},\eta_m)$ at period $m=1$, can be obtained by solving $\mathbf{P}_{m'}$ for some $m'\geq 1$ that uses a data forecast $\{\hat{d}_{\xi}\}_{\xi=1}^{m'}$. Because the forecast is expected to contain errors, the initial targets are expected to be suboptimal but these will be refined as true data is obtained. The proposed scheme also enables the incorporation of forecasts at the beginning of each period to compute the intra-period policies. In particular, given the current guess for the targets $w_{m+1}$, and a forecast $\hat{d}_{m+1}$ over period $m+1$, one can compute the internal control and state policies $\hat{y}_{m+1}=(\hat{x}_{m+1},\hat{u}_{m+1})$ that satisfy the targets $w_{m+1}$ using an intra-period MPC controller. In the linear case, this can be done by solving:
\begin{align} 
& \min_{y\in \mathbb{R}_+^{n_y}} \; \hat{c}^T_{m+1} y \; \textrm{ s.t.} \quad W y = \hat{r}_{m+1} - Tw_{m+1}.
\end{align}
Clearly, the policies $\hat{y}_{m+1}$ are suboptimal because the forecast $\hat{d}_{m+1}$ will deviate from the true realization ${d}_{m+1}$ (this becomes known at the end of period $\xi+1$). Because the cost function $h(\cdot,\cdot)$ is continuous, one can use standard perturbation results to show that the optimality error in the intra-period policy is bounded by the forecast error as $|h(w_{m+1},{d}_{m+1})-h(w_{m+1},\hat{d}_{m+1})|\leq L\|{d}_{m+1}-\hat{d}_{m+1}\|$ for some Lipschitz constant $L\in \mathbb{R}_+$ \cite{shapiropert}. This implies that the quality of the forecast does affect optimality with respect to the intra-period policies. Interestingly, however, the forecast quality does not affect optimality of the targets. 

The full hierarchical scheme is sketched in Figure \ref{fig: summary} and is summarized as follows:
\vspace{-0.1in}
\begin{enumerate}
\item Initialize $m \leftarrow 1$, $\mathbb{V}_{m} \leftarrow \emptyset$, and $w_{m}$.
\item Use forecast $\hat{d}_m$ and targets $w_m$ to obtain intra-period policies $\hat{y}_m$ using MPC controller. 
\item Let system transition from $m\to m+1$.
\item Observe $d_{m}$ and solve $\mathbf{S_m}$ to obtain $\pi(w_{m},d_{m})$. 
\item Update $\mathbb{V}_m\leftarrow \mathbb{V}_{m-1} \cup \{\pi(w_{m},d_{m})\}$. 
\item Obtain $\alpha_\xi^m, \beta_\xi^m,\, \xi =1,\dots,m$ from cut generator.
\item Solve $\mathbf{M}_m$ and obtain updated targets $w_{m+1}$.
\item Shift period time $m\leftarrow m+1$ and return to Step 2.
\end{enumerate}

The SP interpretation allows us to derive metrics and techniques to monitor optimality. We first note that the running cost ${\phi}_m(w_{m})$ evaluated at the current target $w_{m}$ can be evaluated by solving the sequence of subproblems $\mathbf{S}_\xi,\,\xi=1,\dots,m$. The running cost is an upper bound of the optimal running cost (obtained by solving $\mathbf{P}_m$). Moreover, the proposed CP scheme offers the guarantee that the cost $\underline{\phi}_m(w_{m})$ is a lower bound of the running cost ${\phi}_m(w_{m})$, which is an asymptotically exact statistical approximation of $\phi(w_{m})$. The difference between the running cost (upper bound) and the lower bound is known in the SP literature as the optimality gap and is formally defined as $\epsilon_m:=(\phi_m(w_m)-\underline{\phi}_m(w_m))/\phi_m(w_m)$. Here, we refer to $\epsilon_m$ as the current gap. The convergence of the CP scheme guarantees that $\lim_{m\to\infty}\epsilon_m=0$. We note that this gap can be used to measure the quality of the CP approximation but should be used with care when interpreting optimality. In particular, the gap can only be used as a measure of optimality in the limit $m\to \infty$ (once the running cost $\phi_m(w_m)$ is close to the actual cost $\phi_\infty(w_m)$). Consequently, one usually resorts to computing confidence intervals for $\phi_m(w_m)$ by using inference (evaluate the cost at $w_m$ but using different combinations of realizations) \cite{linderoth2006empirical}. Motivated by this, in benchmarking studies, we are also interested in monitoring the overall optimality gap $\bar{\epsilon}_m:=(\phi_\infty(w_m)-\underline{\phi}_m(w_m))/\phi_\infty(w_m)$, provided that $\phi_\infty(w_m)$ can be computed.

\subsection{Extensions to Nonlinear Systems}

In the case of linear systems, the CP scheme can approximate the running cost $\phi_m(w)$ using a finite number of supporting hyperplanes, which keeps the master problem $\mathbf{M}_m$ tractable. The SP representation opens the door to other schemes such as proximal point methods for nonlinear (but convex) problems \cite{bertsekas2011incremental}. Here, the idea is to prevent the accumulation of data over time by summarizing past information in terms of a proximal term. In our context, for instance, the proximal point strategy will result in a problem of the form at stage $m$:
\vspace{-0.1in}
{\setlength{\mathindent}{0pt}
\begin{align}\label{eq:proximal}
\min_{w\in \mathcal{W}} \;  \phantom{ } &\mu \|w-w_{m}\|^2+g(w) +\frac{1}{m'}\sum_{\xi=m'}^m h(w,d_\xi).
\end{align}}
\vspace{-0.3in}

Here, $\mu \|w-w_{m}\|^2$ is a regularization term with parameter $\mu>0$. This term summarizes data before time $m'$ and resembles the arrival cost used in moving horizon estimation \cite{rawlingsmhe}. In the general case of nonconvex problems one can use a statistical approximation $\min_{w\in\mathcal{W}}\; \phi_m(w)$ at every period $m$ by using linear algebra decomposition schemes. In particular, it is well-known that problems with the structure of $\mathbf{P}_m$ give linear algebra systems that enable parallel decomposition \cite{zavalalaird}. Given that the coupling is only in the space of the periodic targets $x_0$, this approach can scale to systems with thousands of states and tens of thousands of periods (realizations). This approach, however, exhibits a fundamental limitation in the number of periods that it can handle (because data is accumulated unboundedly over time). This is, in fact, is also a limitation of statistical approximation schemes for SP. To circumvent this issue, one can use clustering techniques that seek to compress the realization space to maintain a tractable approximation \cite{cao2016clustering}. Such techniques are based on the observation that data realizations tend to be redundant and only a small subset actually impacts the cost. One can quantify the error incurred in the scenario compression by using inference techniques \cite{linderoth2006empirical}. 
\vspace{-0.1in}

%%%%%%%%%%%%%%%%%%%%%%%%%%%%%%%%%%%%%%%%%% 

\section{Computational Experiments} \label{sec:experiments}
\begin{figure}[!htb]
\centering
\includegraphics[width=0.42\textwidth]{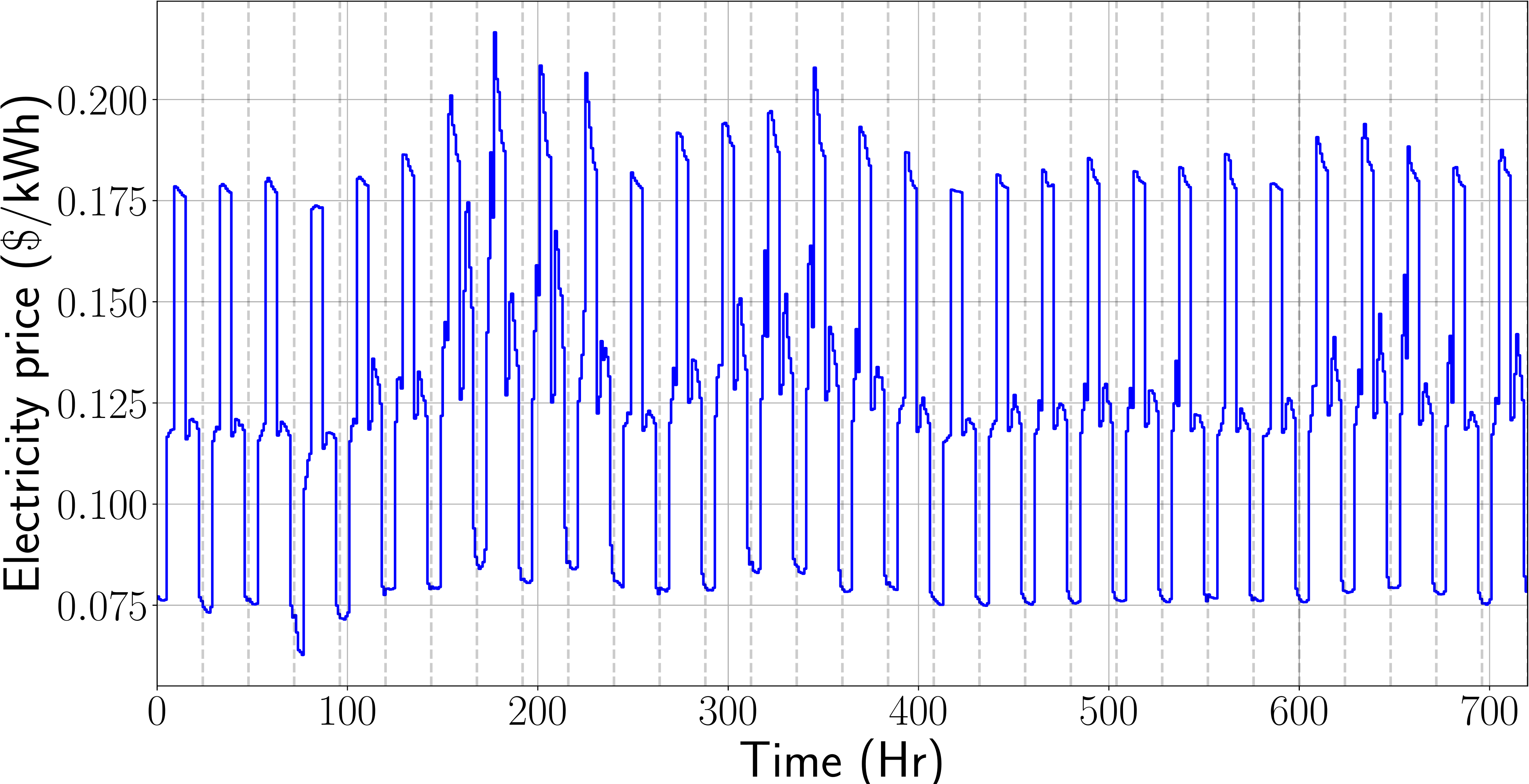} 
\includegraphics[width=0.42\textwidth]{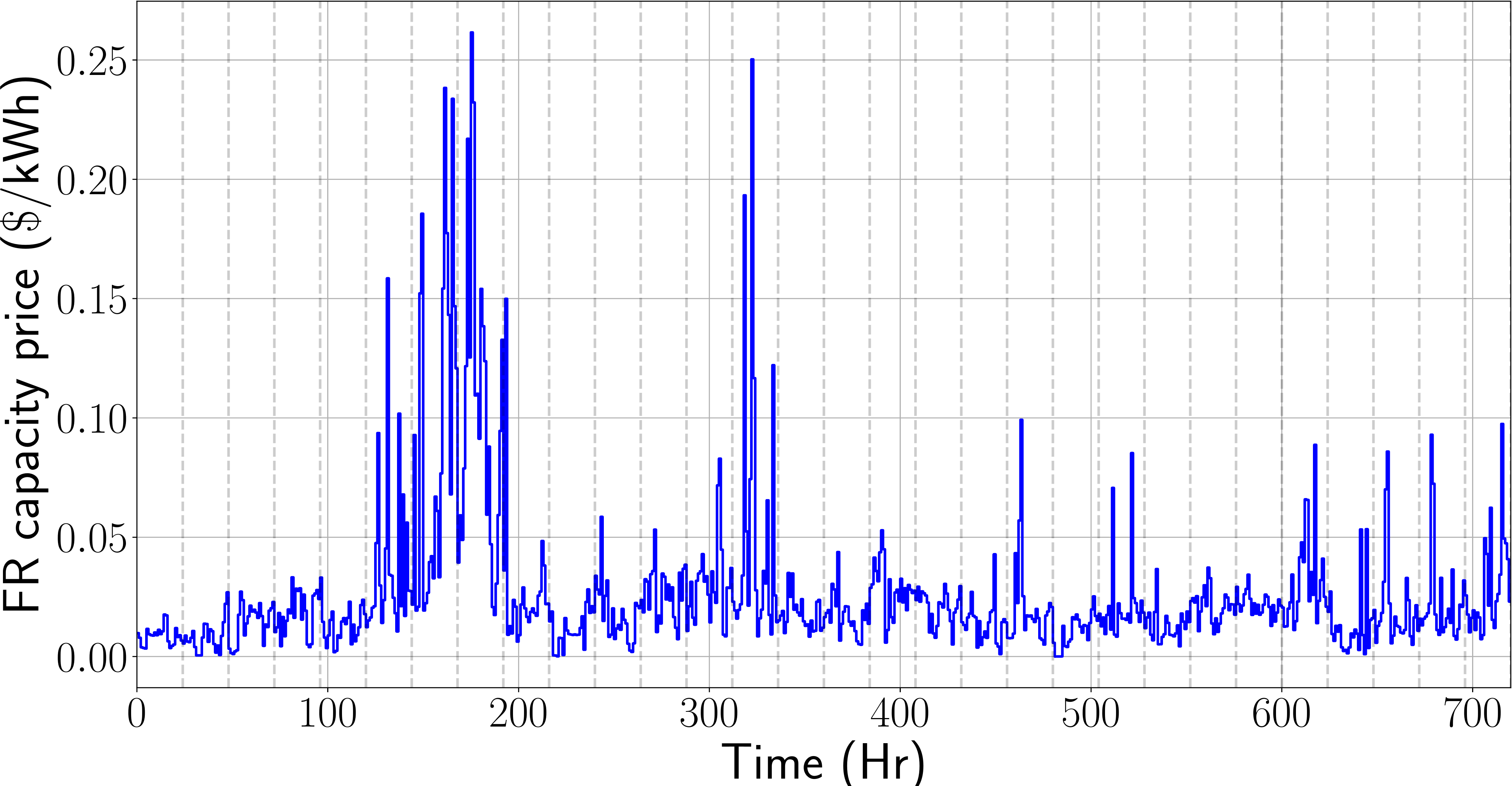}
\includegraphics[width=0.42\textwidth]{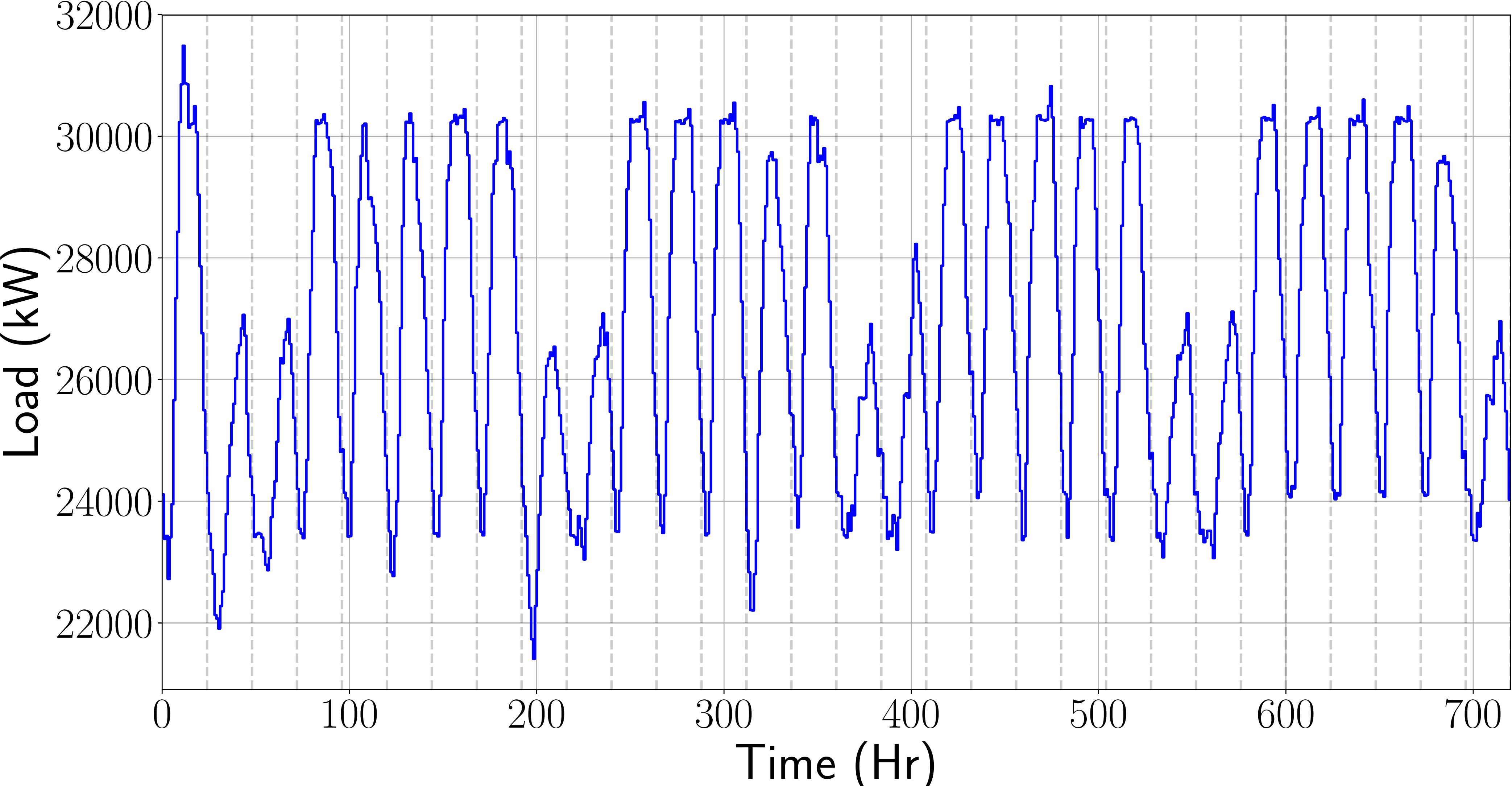} \vspace{-0.1in}
\caption{Energy price (top), FR price (middle), and load (bottom) data.}\label{fig:prices}
\end{figure}

The performance of the proposed scheme is demonstrated using an application in buildings with electricity storage. The goal of the controller is to determine optimal short-term (hourly) participation strategies in frequency regulation (FR) and energy markets that maximize revenue and simultaneously mitigate long-term demand charges. We consider a utility-scale stationary battery with a capacity of 0.5 MWh, rated power of 1 MW, and a ramping limit of 0.5 MW/hr. We use real data for energy and FR prices from PJM Interconnection (shown in Figure \ref{fig:prices}). We also use real load data for a typical university campus (shown also in Figure \ref{fig:prices}). The periodic components in the energy prices and the load profile can be clearly observed, while periodicity of FR prices is not as strong. The MPC problem is formulated using daily periods of $n=24$ hours and we consider an horizon of $m=300$ days (nearly a year). The model parameters include the battery storage capacity ($\overline{E}$ kWh), maximum discharging and charging rates (power) ($\overline{P}, \underline{P}$ in kW, respectively), minimum fraction of battery capacity reserved for frequency regulation ($\rho$ in kWh/kW), and maximum ramping limit ($\overline{\Delta P}$ in kW/h). The random data are the loads ($L_{\xi,t}$ kW), market prices for electricity and FR capacity ($\pi^e_{\xi,t}$ \$/kWh and $\pi^f_{\xi,t}$ \$/kW respectively), demand charge (monthly) ($\pi^D$ in \$/kW), and fraction of FR capacity requested by ISO ($\alpha_{\xi,t}$). The decision variables are net battery discharge rate (power) ($P_{\xi,t}$ in kW), FR capacity provided to ISO ($F_{\xi,t}$ in kW), state-of-charge (SOC) of the battery ($E_{\xi,t}$ in kWh), load requested from utility ($d_{\xi,t}$ in kW) and peak load ($D$ in kW). The formulation minimizes the total cost (negative total revenue), which is given by the demand charge (peak cost) and the revenues collected from the market (time-additive cost). Detailed notation and analysis of this formulation is presented in \cite{zavalaacc}. The problem has the form:
%\begin{subequations}
{\setlength{\mathindent}{0pt}
\begin{align*}
\min \; & \frac{1}{m}\sum\limits_{\xi \in \Xi} \sum\limits_{t \in \mathcal{T}} \left(-\pi^e_{\xi,t}(P_{\xi,t} - \alpha_{\xi,t}F_{\xi,t}) - \pi^f_{\xi,t}F_{\xi,t}\right) + \pi^D D \\
\textrm{s.t. } \; & P_{\xi,t} + F_{\xi,t} \leq \overline{P}, \; t \in \mathcal{T}, \xi \in \Xi\\
& P_{\xi,t} - F_{\xi,t} \geq -\underline{P}, \; t \in \mathcal{T}, \xi \in \Xi  \\
& E_{\xi,t+1} = E_{\xi,t}- P_{\xi,t} + \alpha_{\xi,t}F_{\xi,t}, \; t \in \bar{\mathcal{T}}, \xi \in \Xi\\
& \rho F_{\xi,t} \leq E_{\xi,t}  \leq \overline{E} - \rho F_{\xi,t}, \; t \in \mathcal{T}, \xi \in \Xi \\
& \rho F_{\xi,t} \leq E_{\xi,t+1}  \leq \overline{E} - \rho F_{\xi,t}, \; t \in \bar{\mathcal{T}}, \xi \in \Xi \\
& -\overline{\Delta P}  \leq P_{\xi,t+1} - P_{\xi,t} \leq \overline{\Delta P}, \; t \in \bar{\mathcal{T}}, \xi \in \Xi  \\
& d_{\xi,t} = L_{\xi,t} - P_{\xi,t} + \alpha_{\xi,t} F_{\xi,t}, \; t \in \mathcal{T}, \xi \in \Xi \\
& d_{\xi,t}  \leq D, \; t \in \mathcal{T}, \xi \in \Xi \\
& P_{\xi,t} + F_{\xi,t}  \leq L_{\xi,t}, \; t \in \mathcal{T}, \xi \in \Xi  \\
& E_{\xi,0} = E_{0}, \; \xi \in \Xi \\
& E_{\xi,N_\xi} = E_{0}, \; \xi \in \Xi  \\
& 0  \leq E_{\xi,t}  \leq \overline{E}, \; t \in \mathcal{T}, \xi \in \Xi  \\
& -\underline{P}  \leq P_{\xi,t}  \leq \overline{P},  \;  t \in \mathcal{T}, \xi \in \Xi  \\
& 0  \leq F_{\xi,t}  \leq \overline{P}, \;  t \in \mathcal{T}, \xi \in \Xi
\end{align*}}
\vspace{-.3in}

\begin{figure}[!htb]
\centering
\includegraphics[width=0.40\textwidth]{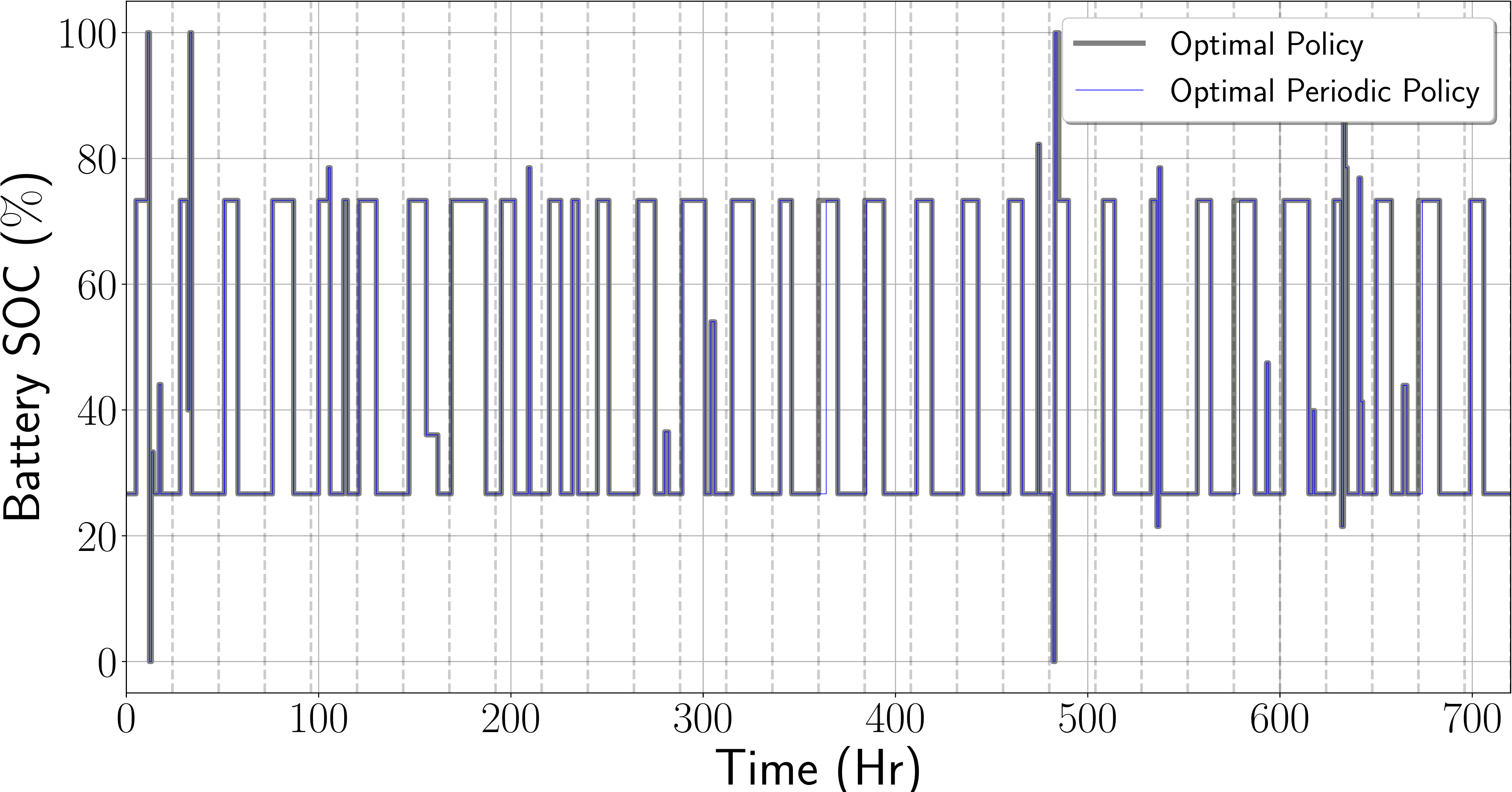} \vspace{-0.1in}
\caption{Optimal and optimal periodic policies.}
\label{fig:perfectsoc}
\end{figure}

\begin{figure}[h!]
\centering
\includegraphics[width=0.47\textwidth]{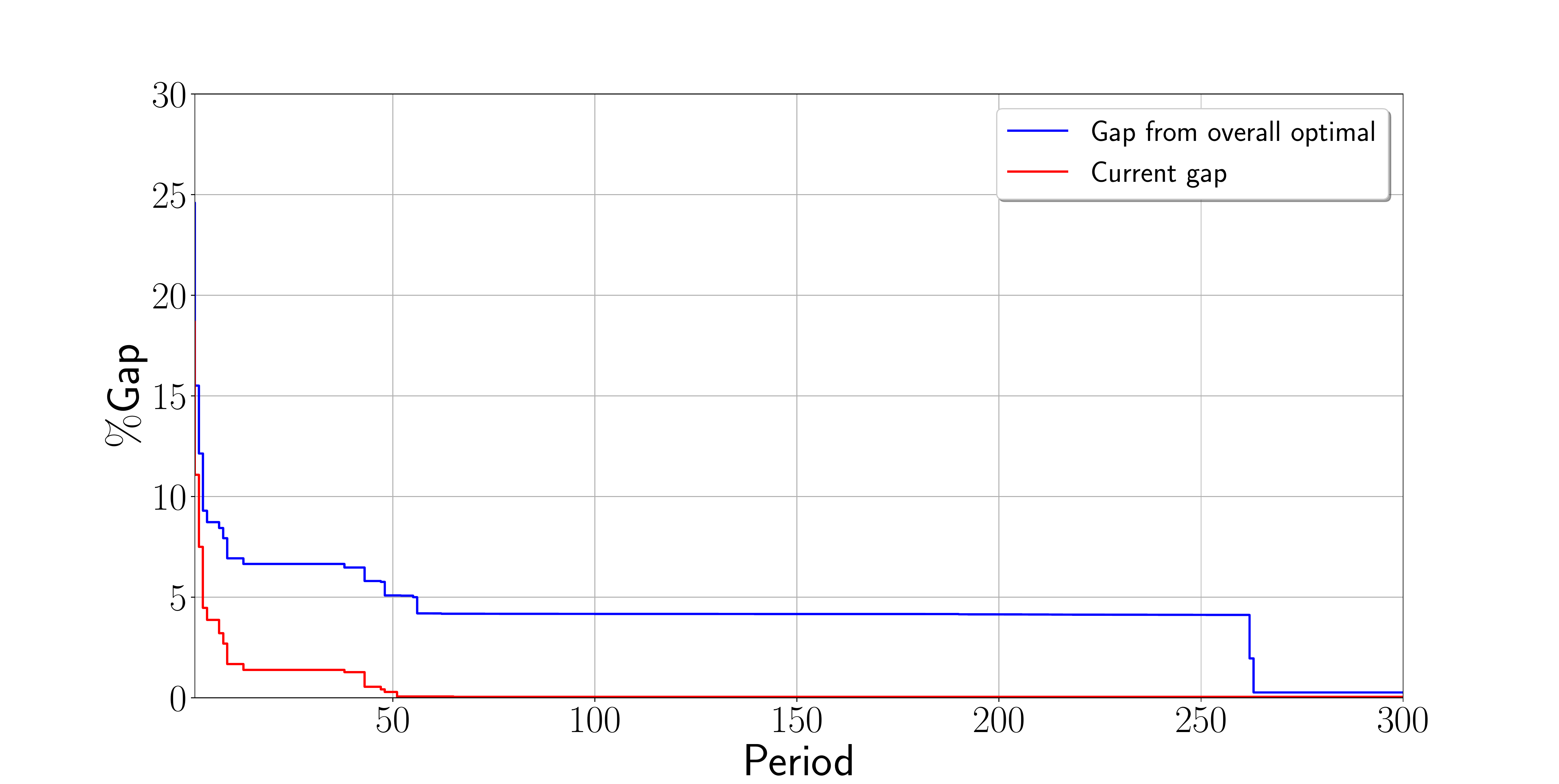} \vspace{-0.1in}
\caption{Evolution of optimality gap.}
\label{fig:gap}
\end{figure}

We first solve the long-horizon MPC problem \eqref{eq:periodicplanning} for a horizon of 30 days assuming perfect knowledge of the data to obtain the optimal policy. We compare this policy against that of a long-horizon MPC formulation that enforces periodicity constraints \eqref{eq:periodic1}. The optimal and periodic policies are presented in Figure \ref{fig:perfectsoc}. Here, we can see that the policies match. In this application, SOC periodicity arises naturally because the battery needs to maintain a minimum SOC level after each period. We ran the proposed retroactive CP scheme for 300 periods to progressively update the targets. The evolution of the current gap $\epsilon_m$ and overall optimality gap $\bar{\epsilon}_m$ is shown in Figure \ref{fig:gap}. We observe that the overall gap eventually vanishes, demonstrating that the CP scheme delivers optimal targets. The overall gap closes to zero close to the end of the horizon, once the peak demand is observed. We also see that the current gap closes to 0\% in about 50 periods and stays there for the rest of the horizon. This illustrates that the cutting planes approximate the cost function well but also that this metric can be misleading. We have also found that the performance of the proposed hierarchical scheme is close to that obtained with the optimal periodic policy (assuming perfect information). In particular, the cost of the optimal policy over the 300 periods is \$136,050 while that of the hierarchical MPC is \$139,978 (a difference of 2.89\%). 

\begin{figure}[h!]
\centering
\includegraphics[width=0.47\textwidth]{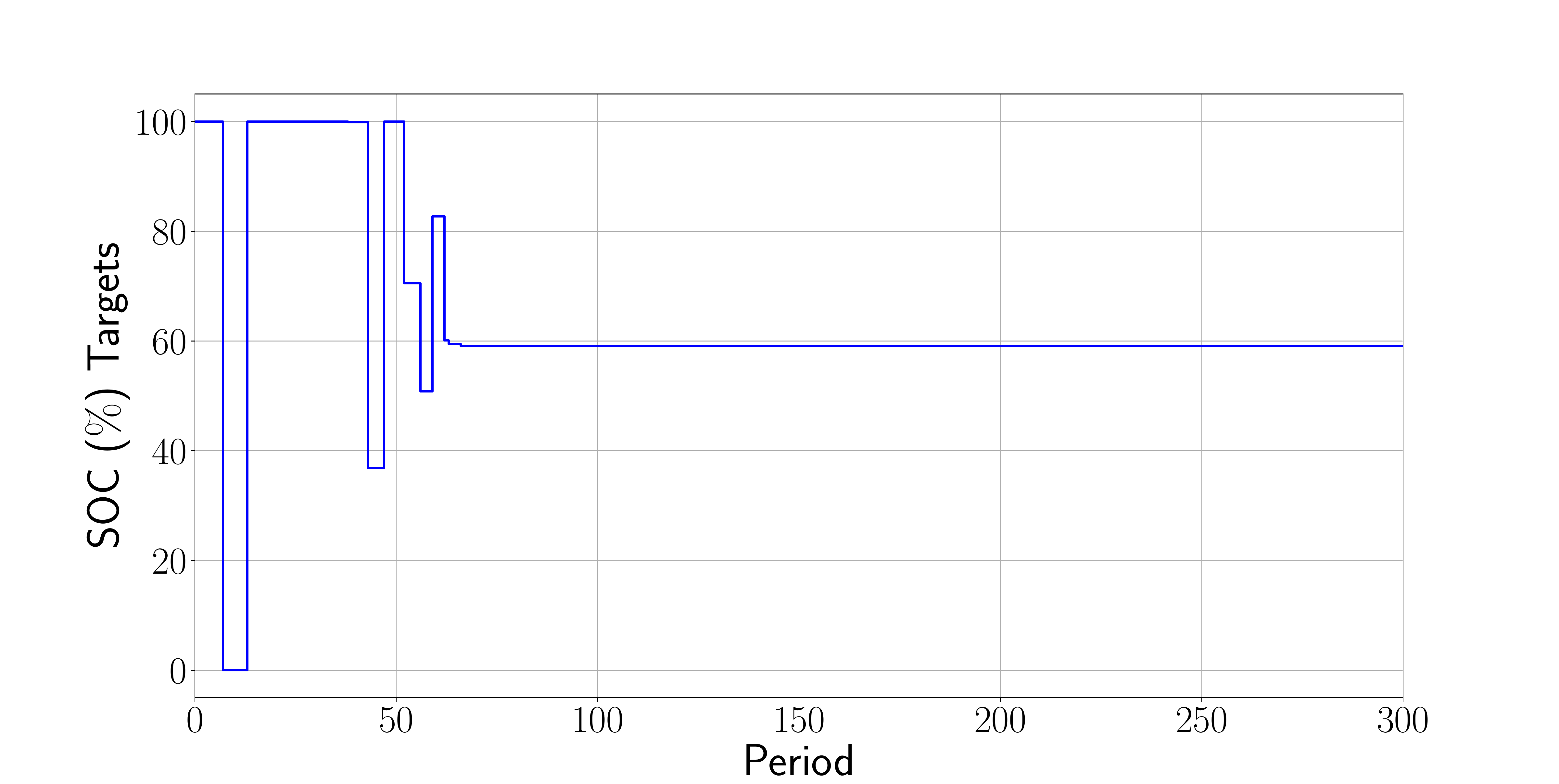} \vspace{-0.1in}
\includegraphics[width=0.47\textwidth]{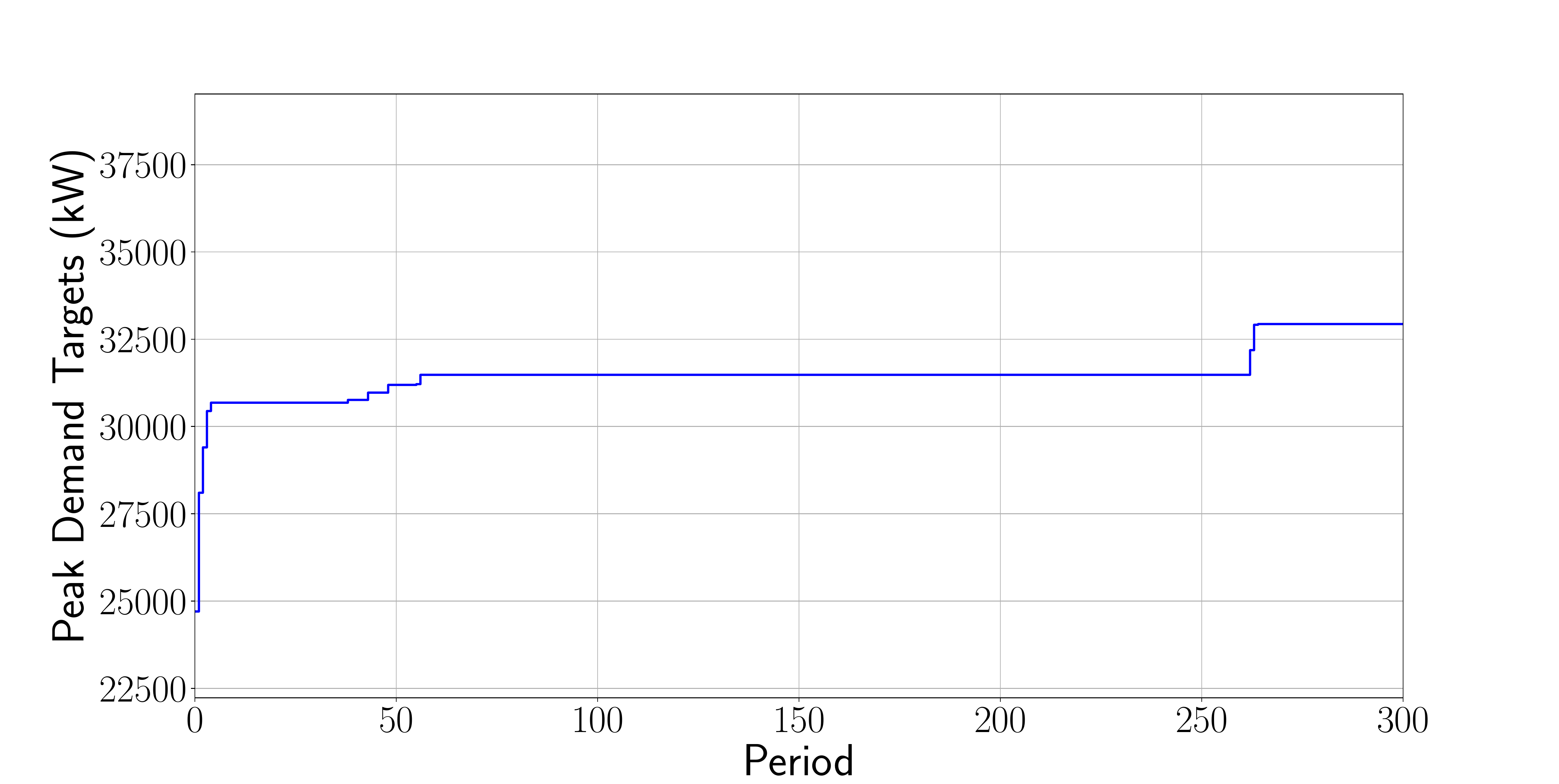}
\caption{Evolution of periodic SOC (top) and peak (bottom) targets obtained with cutting-plane scheme.}
\label{fig:soc}
\end{figure}

\begin{figure}[h!]
\centering
\includegraphics[width=0.47\textwidth]{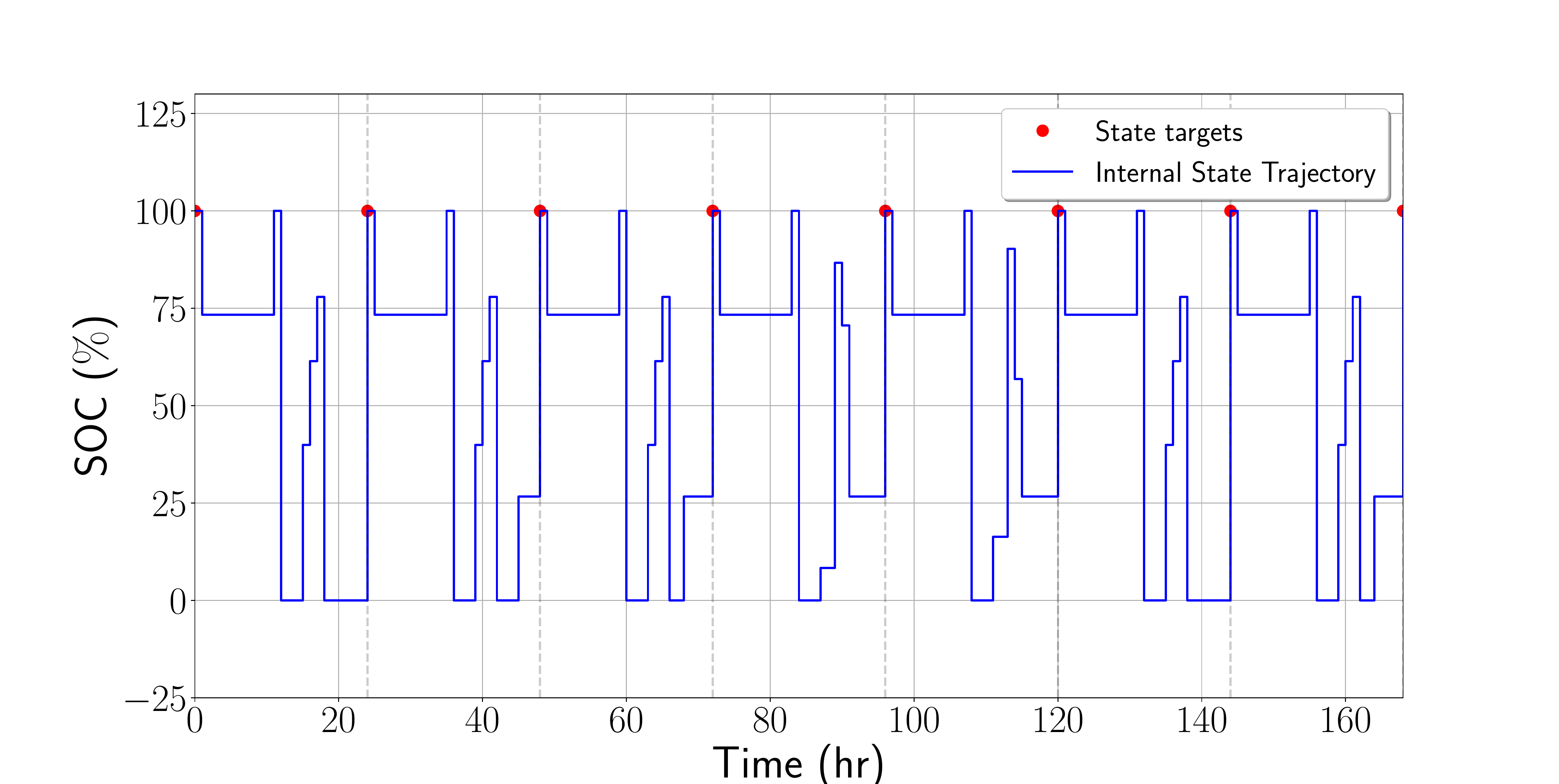} 
\caption{Evolution of periodic SOC targets and intra-period policies obtained with CP scheme (first seven periods).}
\label{fig:traj}
\end{figure}

Figure \ref{fig:soc} shows the evolution of the periodic state and peak targets. We see that the CP scheme adaptively updates the targets as data is accumulated over time. The SOC target settles quickly to the optimal level of 59\%. The peak target requires more periods to settle and this behavior is attributed to the fact that the peak load is observed at period 264. After this period the peak settles at its optimal value of 32,935. Figure \ref{fig:traj} shows the intra-period policies for the short-term MPC controller for the first 7 days (periods) of operation. We see that the controller follows the target of the CP scheme. 


%%%%%%%%%%%%%%%%%%%%%%%%%%%%%%%%%%%%%%%%%%
\vspace{-0.1in}
\section{Conclusions and Future Work}

We showed that stochastic programming provides a framework to design hierarchical MPC schemes for periodic systems. We have shown that, under periodicity, it is possible to compute and refine periodic state targets by solving a retroactive optimization problem that progressively accumulates historical data. The retroactive problem is a statistical approximation of the stochastic program that delivers optimal targets in the long run to guide a short-term MPC controller. This optimality property can be achieved without data forecasts. The SP setting opens the door to a number of potential developments in hierarchical MPC scheme. As part of future work, we are interested in exploring schemes for nonlinear systems and to provide optimality and stability results. Moreover, it is necessary to investigate convergent schemes that prevent accumulation of large amounts of data over time and that factor in forecast information in a more effective manner. 

%%%%%%%%%%%%%%%%%%%%%%%%%%%%%%%%%%%%%%%%%%%%%%%%%%

\bibliographystyle{abbrv}        % Include this if you use bibtex 
\bibliography{hierarchical_mpc_sp}           % and a bib file to produce the 
                                 % bibliography

\end{document}